\documentclass[11pt]{article}
\usepackage{url}

\usepackage{ifpdf}
\ifpdf
    \usepackage[pdftex]{graphicx}   % to include graphics
    \usepackage[pdftex,     % sets up hyperref to use pdftex driver
            plainpages=false,   % allows page i and 1 to exist in the same document
            breaklinks=true,    % link texts can be broken at the end of line
            colorlinks=true,
            pdftitle=MagicSeriesEnum
           ]{hyperref} 
    \usepackage{thumbpdf}
\else
    \usepackage{graphicx}       % to include graphics
    \usepackage{hyperref}       % to simplify the use of \href
\fi

\title{Asymptotic enumeration of magic series}
\author{M. Quist}
\date{\today}

\begin{document}
\maketitle

\begin{abstract}
A magic series is a set of natural numbers that, by virtue of its size, sum, and maximum 
value, could fill a row of a normal magic square.  In this paper,
we derive an exact two-dimensional integral representation for the number of
magic series of order $N$.  By applying the stationary phase approximation,
we develop an expansion in powers of $1/N$ for the number of
magic series and calculate the first few terms.  We find 
excellent agreement between our approximation and the known exact values.
Related results are presented for magic cube and hypercube series, bimagic series,
and trimagic series.
\end{abstract}

\section{Introduction}

A (normal) \emph{magic square} of order $N$ is an arrangement of the integers from $1$ to 
$N^2$ in an $N\times N$ grid such that each row, column, and diagonal has the 
same sum.  The common sum is given by
\begin{equation}
S_1 (N)=\frac{1}{N}\sum_{i=1}^{N^2}i=\frac{1}{2}N (N^2+1).
\end{equation}
A subset of $[N^2]=\{1,2,\ldots,N^2\}$ that could populate a row of a normal magic square
(i.e., an $N$-element subset of $[N^2]$ whose sum is $S_1$) is known as a \emph{magic 
square series}.  Analogously, a (normal)
\emph{magic cube} of order $N$ is an arrangement of the integers from $1$ to $N^3$ into a cube
such that each row, column, pillar, and space diagonal has the same sum;
legal subsets for the rows of a normal magic cube are \emph{magic cube series};
and so on for hypercubes of any order.

The number of magic square series of order $N$ has been computed exactly for $N$ as
large as $1000$~\cite{trump-enum}, and is empirically well-described by H. Bottomley's formula:
\begin{equation}
\label{eq:bottomley}
{\cal N}_{\rm{2}} (N) \approx \left(\frac{1}{\pi} \sqrt{\frac{3}{e}} \right)\cdot \frac{\left(Ne\right)^{N}}{N^3 - \frac{3}{5}N^2 + 
\frac{2}{7}N}.
\end{equation}
While the leading-order behavior can be given a statistical justification, the 
coefficients of the correction terms are purely empirical~\cite{bottomley-eq}.  The goal of the 
present paper is to derive the leading-order behavior and the first two
correction terms analytically, and to present a perturbative framework within which
arbitrarily high-order correction terms can be calculated.

The structure of the paper is as follows.  In Section 2, we describe a 
probabilistic problem that is equivalent, for various parameter values,
to the enumeration of magic series for squares, cubes, and hypercubes.
By applying the central limit theorem to this probabilistic problem, we find 
approximate (Gaussian) formulae for the number of magic series of each type.
In Section 3, we derive an exact two-dimensional integral representation of the desired result,
and show how corrections to the Gaussian formulae can be expressed
as perturbation series in $1/N$.  We compute the first- and second-order
correction terms in these series.  In Section 4, we collect our results
for magic series and compare them to the exact counts that are known.
Finally, in Section 5, we discuss extensions to bimagic and trimagic series,
compute approximations to the number of bimagic square series,
bimagic cube series, and trimagic square series,
and compare these results to known exact counts.

\section{Probabilistic approach}
\subsection{Formulation}
For each $i \ge 1$, let $(X_i, Y_i)$ be a random vector that is equal to $(1, i)$ with probability
$\beta$ and equal to $(0, 0)$ otherwise.
Then the partial sums $(A_m, B_m) = \sum_{i=1}^{m} (X_i, Y_i)$ have the following
interpretation: the probability that $(A_m, B_m)$ is equal to $(A, B)$ is the probability of 
selecting $A$ numbers and arriving at a total of $B$, selecting each value in $\{1,2,\ldots,m\}$ 
for inclusion with probability $\beta$.
When we convert this probability to a \emph{count}, the relevance to enumerating magic
series is clear.  Specifically,
\begin{eqnarray}
\label{statform}
p (A, B) &\equiv& P[ (A_m, B_m) = (A,B) ] \nonumber \\
&=& \beta^{A} (1-\beta)^{m-A} {\cal N} (m, A, B),
\end{eqnarray}
where ${\cal N}(m, A, B)$ is the number of subsets of $\{1, 2, \ldots, m\}$ that have
$A$ elements and sum to $B$.  For appropriately chosen arguments, the 
value of $\cal{N}$ is equal to the number of magic series for squares, cubes, or 
hypercubes.

For magic series, we will be interested in cases where $A_m$ and $B_m$
take on their mean values, which are
\begin{eqnarray}
\mu_x &=& E[A_m] = \beta \sum_{i=1}^{m} 1 = \beta m \\
\mu_y &=& E[B_m] = \beta \sum_{i=1}^{m} i = \frac{1}{2} \beta m (m + 1).
\end{eqnarray}
Define the ``central probability" to be
\begin{equation}
P^{*}(m, \beta) = p(\mu_x, \mu_y) = p\left(\beta m, \frac{1}{2}\beta 
m(m+1)\right).
\end{equation}
The corresponding count is given by
\begin{eqnarray}
\label{prefactor}
{\cal N}^{*}(m, \beta) &=& {\cal N} \left(m, \beta m, \frac{1}{2}\beta m(m+1)\right)  
\nonumber
\\ &=& \beta^{-\beta m}(1-\beta)^{-(1-\beta)m} P^{*}(m, \beta).
\end{eqnarray}
The prefactor in Eq.~(\ref{prefactor})
is exact.  We will now develop an approximation scheme for
the central probability $P^{*}(m, \beta)$, beginning with the Gaussian approximation.

\subsection{Gaussian approximation}

As $m$ becomes large, the central limit theorem (the Lyapunov central limit
theorem, specifically, since the summed variables are independent but not identically
distributed) implies that the distribution of $(A_m, B_m)$ approaches a two-dimensional Gaussian.
The center of the Gaussian is $(\mu_x, \mu_y)$.
The components of its covariance matrix are
\begin{eqnarray}
\Sigma_{xx} &=& E[A_m^2] - E[A_m]^2 \nonumber
\\ &=& \beta (1-\beta) m
\\
\Sigma_{xy} &=& E[A_m B_m] - E[A_m] E[B_m]  \nonumber
\\ &=& \frac{1}{2} \beta (1-\beta) m (m + 1)
\\
\Sigma_{yy} &=& E[B_m^2] - E[B_m]^2  \nonumber
\\ &=& \frac{1}{6} \beta (1-\beta) m (m + 1) (2m + 1).
\end{eqnarray}
The limiting distribution of $(A_m, B_m)$ is therefore
\begin{equation}
p_{g}(A,B) = \frac{1}{2\pi\sqrt{\det\Sigma}}\exp(\cdots) = \frac{1}{\pi \beta(1-\beta) m} \sqrt{\frac{3}{m^2-1}} 
\exp(\cdots),
\end{equation}
where the exponent is quadratic in the distance from $(A, B)$ to the center of the 
distribution, and we have used the following for the square root of the determinant
of the covariance matrix:
\begin{eqnarray}
\sqrt{\det \Sigma} &=& \beta (1-\beta) \sqrt{ \frac{1}{6} m^2 (m+1) (2m+1) - \frac{1}{4} m^2 (m+1)^2 } \nonumber
\\ &=& \frac{1}{2} \beta (1-\beta) m \sqrt{\frac{m^2 - 1}{3}}.
\end{eqnarray}
In the Gaussian approximation, then, the central probability is
\begin{equation}
P^{*}_{g}(m,\beta) = \frac{1}{\pi \beta(1-\beta) m} \sqrt{\frac{3}{m^2-1}}.
\end{equation}
We can obtain results for various interesting cases by
specifying $m$ and $\beta$.

\subsection{Magic square series}

The number of magic square series of order $N$ is given by ${\cal N}_2(N) = {\cal N}(N^2, N, 
S_1)$, where $S_1=\frac{1}{2}N(N^2+1)$.
We see that if we choose $m=N^2$ and $\beta = 1/N$, then the center of the 
Gaussian is located at $(\mu_x, \mu_y) = (N, S_1)$.  In this case, the prefactor in Eq.~({\ref{prefactor}}) is
\begin{eqnarray}
\label{prefactor-squares}
\beta^{-\beta m}(1-\beta)^{-(1-\beta)m} &=& N^{N}\left(1-\frac{1}{N}\right)^{-(N^2-N)}  \nonumber
\\ &=& N^N \exp\left[-(N^2-N)\ln\left(1-\frac{1}{N}\right)\right] \nonumber
\\ &=& N^N \exp\left[N - \frac{1}{2} - \frac{1}{6N} - \frac{1}{12N^2} + \cdots\right] 
\nonumber
\\ &=& \frac{\left(Ne\right)^{N}}{\sqrt{e}} \exp\left[-\frac{1}{6N}-\frac{1}{12N^2}+\cdots\right] 
\nonumber
\\ &=& \frac{\left(Ne\right)^{N}}{\sqrt{e}} \left(1 - \frac{1}{6N} - \frac{5}{72N^2} + 
\cdots\right),
\end{eqnarray}
where the series is readily extended to any order.  The Gaussian approximation for the number of 
magic square series of order $N$ is therefore
\begin{eqnarray}
{\cal N}_{2,g}(N) &=&
\frac{\left(Ne\right)^{N}}{\sqrt{e}} \left(1 - \frac{1}{6N} - \frac{5}{72N^2} + 
\cdots\right)\times P_{g}^{*}(N^2, 1/N) \nonumber \\
&=&
\frac{\left(Ne\right)^{N}}{\sqrt{e}} \left(1 - \frac{1}{6N} - \frac{5}{72N^2} + 
\cdots\right)\times \frac{1}{\pi (N-1)} \sqrt{\frac{3}{N^4-1}} \nonumber
\\ &=& \left(\frac{1}{\pi} \sqrt{\frac{3}{e}} \right)\cdot \frac{\left(Ne\right)^{N}}{N^3} \left(1 + \frac{5}{6N} + \frac{55}{72N^2} + \cdots\right).
\end{eqnarray}

When we compare this estimate to the actual counts for large $N$,
we find that it is
correct to leading order, but the higher-order correction coefficients are wrong.  (For instance, the coefficient of $1/N$
appears to be $3/5$, not $5/6$.)  This should not be terribly surprising.  First, although the CLT
guarantees a Gaussian distribution in the limit $m\rightarrow\infty$, there are 
finite-$m$ corrections that affect the higher-order terms.  Second, as the careful reader will have noted,
the CLT no longer applies when $\beta$ is a function of $m$.  We can only guarantee convergence 
to a normal distribution when $\beta$ is fixed.  To address both these problems, 
we will need to systematically include non-Gaussian corrections, as we will do in Section 3.

\subsection{Magic cube series}

The number of magic cube series of order $N$ is given by ${\cal N}_3(N) = {\cal N}(N^3, N, 
T_1)$, where
\begin{equation}
T_1(N) = \frac{1}{N^2}\sum_{i=1}^{N^3}i = \frac{1}{2}N(N^3+1).
\end{equation}
If we choose $m=N^3$ and $\beta=1/N^2$, then the center of the Gaussian is 
located at $(\mu_x, \mu_y) = (N, T_1)$.  The prefactor in this case is
\begin{eqnarray}
\label{prefactor-cubes}
\beta^{-\beta m}(1-\beta)^{-(1-\beta)m} &=& N^{2N}\left(1-\frac{1}{N^2}\right)^{-(N^3-N)}  
\nonumber
\\ &=& N^{2N} \exp\left[-(N^3-N)\ln\left(1-\frac{1}{N^2}\right)\right] \nonumber
\\ &=& N^{2N} \exp\left[N - \frac{1}{2N} - \frac{1}{6N^3} + \cdots\right] 
\nonumber
\\ &=& \left(N^2 e\right)^{N} \exp\left[-\frac{1}{2N}-\frac{1}{6N^3}+\cdots\right] 
\nonumber
\\ &=& \left(N^2 e\right)^{N} \left(1 - \frac{1}{2N} + \frac{1}{8N^2} + 
\cdots\right),
\end{eqnarray}
and the Gaussian approximation for the number of magic cube series 
of order $N$ is
\begin{eqnarray}
{\cal N}_{3,g}(N) &=& \left(N^2 e\right)^{N}\left(1 - \frac{1}{2N} + \frac{1}{8N^2} + 
\cdots\right) \times P_{g}^{*}(N^3, 1/N^2) \nonumber
\\ &=&
\left(N^2 e\right)^{N}\left(1 - \frac{1}{2N} + \frac{1}{8N^2} + 
\cdots\right) \times \frac{1}{\pi \left(N-\frac{1}{N}\right)} \sqrt{\frac{3}{N^6-1}} 
\nonumber
\\ &=& \left(\frac{\sqrt{3}}{\pi} \right) \cdot \frac{\left(N^2 e\right)^{N}}{N^4} 
\left(1 - \frac{1}{2N} + \frac{9}{8N^2} + \cdots\right).
\end{eqnarray}
As with the magic square series, a comparison to the exact counts shows that the leading order here is 
correct, while the higher-order terms are wrong.

\subsection{Magic series for $\alpha$-hypercubes}

Generalizing the above by taking $m=N^{\alpha}$ and $\beta=1/N^{\alpha-1}$, we obtain an estimate for
the number of magic series for $\alpha$-hypercubes.  Assuming $\alpha > 3$, the prefactor becomes
\begin{eqnarray}
\beta^{-\beta m}(1-\beta)^{-(1-\beta)m} &=& N^{(\alpha-1) N}\left(1-\frac{1}{N^{\alpha-1}}\right)^{-(N^\alpha-N)} 
\nonumber
\\ &=& N^{(\alpha-1)N} \exp\left[-(N^\alpha-N)\ln\left(1-\frac{1}{N^{\alpha-1}}\right)\right] 
\nonumber
\\ &=& N^{(\alpha-1)N} \exp\left[N - \frac{1}{2N^{\alpha-2}} + \cdots\right] 
\nonumber
\\ &=& \left(N^{\alpha-1} e\right)^{N} \exp\left[-\frac{1}{2N^{\alpha-2}} +\cdots\right] 
\nonumber
\\ &=& \left(N^{\alpha-1} e\right)^{N} \left(1 - \frac{1}{2N^{\alpha-2}} + \cdots\right).
\end{eqnarray}
The Gaussian approximation for the number of magic series for $\alpha$-hypercubes of order $N$
($\alpha > 3$) is
\begin{eqnarray}
{\cal N}_{\alpha>3,g}(N) &=& \left(N^{\alpha-1} e\right)^{N}\left(1 - \frac{1}{2N^{\alpha-2}} + 
\cdots\right) \times P_{g}^{*}(N^\alpha, 1/N^{\alpha-1}) \nonumber
\\ &=&
\left(N^{\alpha-1} e\right)^{N}\left(1 - \frac{1}{2N^{\alpha-2}} +
\cdots\right) \times \frac{1}{\pi \left(N-\frac{1}{N^{\alpha-1}}\right)} \sqrt{\frac{3}{N^{2\alpha}-1}} 
\nonumber
\\ &=& \left(\frac{\sqrt{3}}{\pi} \right) \cdot \frac{\left(N^{\alpha-1} e\right)^{N}}{N^{\alpha + 1}} 
\left(1 - \frac{1}{2N^{\alpha-2}} + \cdots\right).
\end{eqnarray}
Again, only the leading-order behavior is correct.

\section{Integral representation}

\subsection{Formulation}
Now consider evaluating the probability distribution $p(x, y)$ without
appealing to the central limit theorem.  The probability distribution is a 
convolution of $m$ simple distributions; in Fourier space, this convolution 
becomes a product.  Defining $\tilde{p}(k_x, k_y) = \sum_{x,y}\exp(i{\bf k}\cdot{\bf x}) 
p(x, y)$, we have
\begin{equation}
\tilde{p}(k_x, k_y) = \prod_{j=1}^{m} \tilde{p_j}(k_x, k_y),
\end{equation}
where
\begin{eqnarray}
\tilde{p_j}(k_x, k_y) &=& (1-\beta) + \beta e^{i z_j} \nonumber
\\ &=& \exp\left[\ln \left(1 - \beta(1-e^{iz_j})\right) \right],
\end{eqnarray}
and $z_j=(1, j)\cdot {\bf k}= k_x + jk_y$.  Taking the product over $j$
gives
\begin{equation}
\tilde{p}(k_x, k_y) = \exp\left[\sum_{j=1}^{m}\ln \left(1 - 
\beta(1-e^{iz_j})\right)\right].
\end{equation}
The central probability is then given by the inverse Fourier
transform:
\begin{eqnarray}
\label{exact-integral-rep}
P^{*}(m, \beta) &=& p(\mu_x, \mu_y) \nonumber \\
&=& \int_{-\pi}^{\pi}\int_{-\pi}^{\pi}\frac{dk_x dk_y}{(2\pi)^2}\exp\left(-i\mu_x k_x-i\mu_y k_y\right)
\tilde{p}(k_x, k_y) \nonumber \\
&=&
\int\frac{d^2 k}{(2\pi)^2}\exp\left\{\sum_{j=1}^{m}\left[-i\beta z_j + \ln\left(1-\beta(1-e^{iz_j})\right)\right]\right\}.
\end{eqnarray}
This expression, while exact, is not of very much use as it stands, and
will need to be approximated.

\subsection{Stationary phase approximation}

First, expanding in $\beta$, we rewrite the exponent as
\begin{equation}
\beta\sum_{j=1}^{m}(e^{iz_j}-1-iz_j) - 
\sum_{k=2}^{\infty}\frac{\beta^k}{k}\sum_{j=1}^{m}(1-e^{iz_j})^k.
\end{equation}
This expansion is convergent for small $\beta$.  We now propose that the primary contribution to
the integral in Eq.~(\ref{exact-integral-rep}) comes from the region around the origin, where the
real part of the exponent is zero and its imaginary part is stationary with 
respect to $k$.  
This leads us to employ the \emph{stationary phase approximation}, in which the integral of a rapidly oscillating function is approximated
by the local integral(s) at the point(s) of stationary phase.  As the oscillations in
phase become increasingly rapid (which we expect to occur here as $m\rightarrow \infty$),
contributions from all other points are exponentially suppressed.  (In this case, the amplitude of
the integrand is also exponentially suppressed away from the origin.)

Collecting like powers of $z$, we have
\begin{equation}
-\frac{1}{2}\beta(1-\beta)\sum_{j=1}^{m}z_j^2 + \beta\sum_{k=3}^{\infty} \frac{f_k(\beta) \sum_{j=1}^{m}\left(iz_j\right)^k}{k!},
\end{equation}
in the exponent, where each ``vertex correction" function $f_k(\beta)$ is equal to $1$ plus
corrections of order $\beta$.  To be more precise, the first few vertex corrections are
\begin{eqnarray}
f_3(\beta) &=& 1 - 3\beta + 2\beta^2 \\
f_4(\beta) &=& 1 - 7\beta + 12\beta^2 - 6\beta^3.
\end{eqnarray}
Summing over $j$ then gives the exponent as
\begin{eqnarray}
- \frac{\Sigma_{ab}}{2}k_{a}k_{b} - \frac{i \beta f_3(\beta) T_{abc}}{6} k_{a} k_{b} k_{c} 
+ \frac{\beta f_4(\beta) S_{abcd}}{24} k_{a}k_{b}k_{c}k_{d} + \ldots,  
\end{eqnarray}
where repeated indices are summed over, $\Sigma$ is exactly the same as before, and
\begin{eqnarray}
T_{xxx} &=& \sum_{j=1}^{m} 1 = m  \\
T_{xxy} &=& \sum_{j=1}^{m} j =\frac{1}{2} m(m+1) \\
T_{xyy} &=& \sum_{j=1}^{m} j^2 =\frac{1}{3} m(m+1)\left(m+\frac{1}{2}\right) \\
T_{yyy} &=& \sum_{j=1}^{m} j^3 =\frac{1}{4} m^2(m+1)^2.
\end{eqnarray}
All higher-order vertices follow the same pattern: the value of $V_{abc\ldots}$, 
where $V$ is any vertex and $abc\ldots$ are its indices,
is equal to $\sum_{j=1}^{m} j^{n_y(abc\ldots)}$, where $n_y$ counts the number 
of $y$'s among the indices.  To leading order in $m$,
\begin{equation}
  V_{abc\ldots} = \frac{m^{1 + n_ y(abc\ldots)}}{1 + n_y(abc\ldots)}.
\end{equation}
The central probability is now given by
\begin{eqnarray}
P^{*}(m, \beta) &=&
\int\frac{d^2 k}{(2\pi)^2}\exp\Big[- \frac{\Sigma_{ab}}{2}k_{a}k_{b} - \frac{i\beta f_3(\beta) V_{abc}}{6} k_{a} k_{b} 
k_{c} \nonumber
\\ &&\qquad\qquad+ \frac{\beta f_4(\beta) V_{abcd}}{24} k_{a}k_{b}k_{c}k_{d} + \ldots\Big].
\end{eqnarray}
Using rescaled variables $\tilde{k}_x=\sqrt{m\beta(1-\beta)} k_x$
and $\tilde{k}_y=m\sqrt{m\beta(1-\beta)} k_y$, we remove all $\beta$-dependence 
and the leading-order dependence on $m$ from the vertex terms, giving
\begin{eqnarray}
&& P^{*}(m, \beta) = \frac{1}{\beta(1-\beta)m^2} \int\frac{d^2 \tilde{k}}{(2\pi)^2}\exp\Big[- \frac{\tilde{\Sigma}_{ab}}{2}\tilde{k}_{a}\tilde{k}_{b}
\nonumber \\ && \qquad\qquad\qquad - \frac{i \tilde{f}_3(\beta) \tilde{V}_{abc}}{6 \sqrt{\beta m}}  \tilde{k}_{a} \tilde{k}_{b} 
\tilde{k}_{c} +
 \frac{\tilde{f}_4(\beta)\tilde{V}_{abcd}}{24 \beta m} \tilde{k}_{a}\tilde{k}_{b}\tilde{k}_{c}\tilde{k}_{d} 
+ \ldots\Big].
\end{eqnarray}
Here
\begin{equation}
\tilde{\Sigma} = \left(
\begin{array}{cc}
  1 & \frac{1}{2}(1+\frac{1}{m}) \\
  \frac{1}{2}(1+\frac{1}{m}) & \frac{1}{3} (1+\frac{1}{m})(1+\frac{1}{2m})
\end{array}
\right)
\end{equation}
exactly, and
\begin{equation}
\tilde{V}_{abc\ldots} = \frac{1}{1 + n_y(abc\ldots)}
\end{equation}
up to corrections of order $1/m$, and $\tilde{f}_k(\beta) = 
(1-\beta)^{-k/2}f_k(\beta)$.  The first two rescaled vertex corrections are
\begin{eqnarray}
\tilde{f}_3(\beta) &=& 1 - \frac{3}{2}\beta + O(\beta^2) \\
\tilde{f}_4(\beta) &=& 1 - 5\beta + O(\beta^2).
\end{eqnarray}
For simplicity we will drop the tildes on $k$ and work with the rescaled variables from this point
forward.  At this point we have arrived at our main theoretical result:
\begin{eqnarray}
\label{perturb-series}
&&P^{*}(m, \beta) = P_{g}^{*}(m, \beta) \times \nonumber \\ &&\qquad\left\langle \exp\left(-\frac{i\tilde{f}_3(\beta)\tilde{V}_{abc}}{6\sqrt{\beta m}}k_a k_b k_c + \frac{\tilde{f}_4(\beta)\tilde{V}_{abcd}}{24\beta m}k_a k_b k_c k_d + \ldots\right) \right\rangle 
\end{eqnarray}
where the expectation value is evaluated with respect to the normalized Gaussian 
distribution with covariance matrix $\tilde{\Sigma}^{-1}$.  Because each additional power of $k$ (that is,
each new vertex) comes with a multiplier of $(\beta m)^{-1/2}$, expanding the exponential gives
us an asymptotic series, in powers of $1/(\beta m)$, for the multiplicative correction
to $P_{g}^{*}(m, \beta)$.  Odd powers of $k$ evaluate to zero by symmetry, so 
only integer powers of $1/(\beta m)$ appear.

\subsection{First-order correction}

To first order in $1/(\beta m)$, we have
\begin{equation}
  P^{*}(m, \beta) \approx P_g^{*}(m, \beta)\left(1 + \frac{K_1}{\beta m}\right),
\end{equation}
where
\begin{equation}
K_1 = -\frac{\tilde{V}_{abc} \tilde{V}_{def}}{72}\langle k_a k_b k_c k_d k_e k_f \rangle + \frac{\tilde{V}_{abcd}}{24}\langle k_a k_b k_c k_d 
\rangle.
\end{equation}
Note that the cubic and quartic (in $k$) corrections both contribute to this order.
The usual diagrammatic approach applies: the ``propagator" is
\begin{equation}
  \langle k_a k_b \rangle = \Pi_{ab},
\end{equation}
where
\begin{equation}
  \Pi \equiv {\tilde{\Sigma}}^{-1} = \left(
  \begin{array}{cc}
    4 & -6 \\ -6 & 12
  \end{array}\right)
\end{equation}
to leading order
in $1/m$, and expectation values of higher powers of $k$ can be evaluated by combining pairs
of $k$'s in every possible way, e.g.,
\begin{equation}
\langle k_a k_b k_c k_d \rangle = \Pi_{ab} \Pi_{cd} + \Pi_{ac} \Pi_{bd} + \Pi_{ad} 
\Pi_{bc}.
\end{equation}
Because $\tilde{V}$ is symmetric on its indices, this gives
\begin{equation}
\label{eq:k1}
K_1 = -\frac{\tilde{V}_{abc} \tilde{V}_{def}}{8} \Pi_{ab}\Pi_{cd}\Pi_{ef} -\frac{\tilde{V}_{abc} \tilde{V}_{def}}{12}\Pi_{ad}\Pi_{be}\Pi_{cf} + \frac{\tilde{V}_{abcd}}{8} \Pi_{ab} \Pi_{cd} .
\end{equation}
In terms of diagrams, we are summing over all distinct diagrams
for which the sum of the vertex degrees is equal to twice the number
of vertices plus two.  There are three
such diagrams, corresponding to the three terms in Eq.~(\ref{eq:k1}).
A straightforward summation over the indices yields the correction 
coefficient:
\begin{equation}
K_1 = -\frac{1}{8}\cdot 4 -\frac{1}{12}\cdot 4 + \frac{1}{8}\cdot\frac{24}{5} = -\frac{7}{30}.
\end{equation}

\subsection{Higher-order corrections}

The coefficient $K_1$ derived in the previous section is
\emph{universal}, in that it describes the first non-Gaussian correction to $P^{*}$ 
whenever $\beta\rightarrow 0$ and $(\beta m) \rightarrow \infty$.  In 
particular, the same value of $K_1$ applies in counting magic 
series for hypercubes of any order.  In each case, the Gaussian approximation 
derived in Section 2 must be multiplied by $1 + K_1/(\beta m) = 1 - 7/(30N)$ to 
yield the correct expression to order $1/N$.

Higher-order corrections are not universal in the same way, because there are 
correction terms in powers of $1/m$ (from approximations made to the propagator and 
the vertex), powers of $1/(\beta m)$ (from higher-order diagrams), and 
powers of $\beta$ (from the vertex corrections).  Depending on the
dimensionality of the problem, these correction terms scale differently.  For magic
square series, $\beta = 1/(\beta m)=1/N$, and so we need both the next-order terms in $\beta$ 
and the next round of diagrams to get the coefficient of $1/N^2$ in the large-$N$ expansion.
For magic cube and hypercube series, on the other hand, $\beta \ll 1/N$, and so
only the next round of diagrams are required for order $1/N^2$.  (The approximations made to the
propagator and the vertex produce still finer corrections: they enter the calculation at order $1/N^3$
for magic series for squares, and at even higher orders for cubes and hypercubes.)

So, to correctly determine the next order in $1/N$ we need at most two types of additional correction.  Omitting the 
details, we find
\begin{equation}
P^{*}(m, \beta) \approx P^{*}_{g}(m, \beta) \left(1 + \frac{K_1}{\beta m} + \frac{K_2}{m} + \frac{K_3}{(\beta m)^2}\right),
\end{equation}
where
\begin{eqnarray}
K_2 &=& \frac{3\tilde{V}_{abc} \tilde{V}_{def}}{8} \Pi_{ab}\Pi_{cd}\Pi_{ef} +\frac{\tilde{V}_{abc} \tilde{V}_{def}}{4}\Pi_{ad}\Pi_{be}\Pi_{cf} - \frac{5\tilde{V}_{abcd}}{8} \Pi_{ab} \Pi_{cd}  
\nonumber \\
&=& \frac{3}{8}\cdot 4 + \frac{1}{4}\cdot 4 - \frac{5}{8}\cdot\frac{24}{5} 
\nonumber \\
&=& -\frac{1}{2}
\end{eqnarray}
arises from $\tilde{f}_k(\beta)$-corrections, and includes only the three diagrams 
already considered; and
\begin{eqnarray}
K_3 &=&
\frac{\tilde{V}_{abc}\tilde{V}_{def}\tilde{V}_{ghi}\tilde{V}_{jkl}}{31104}\langle k_{abcdefghijkl} \rangle
- \frac{\tilde{V}_{abc}\tilde{V}_{def}\tilde{V}_{ghij}}{1728}\langle k_{abcdefghij}\rangle
\nonumber \\ &&
+ \frac{\tilde{V}_{abc} \tilde{V}_{defgh}}{720} \langle k_{abcdefgh}\rangle
+ \frac{\tilde{V}_{abcd}\tilde{V}_{efgh}}{1152}\langle k_{abcdefgh}\rangle
- \frac{\tilde{V}_{abcdef}}{720}\langle k_{abcdef}\rangle
\nonumber \\  &=& \frac{205}{72} - \frac{31}{6} + \frac{22}{15} + \frac{17}{15} - \frac{29}{105}
\nonumber \\  &=& \frac{11}{2520}
\end{eqnarray}
includes a large number of higher-order diagrams.  The summation was performed with the help
of a computer.

\section{Magic series results}

After multiplying the Gaussian results from Section 2 by the correction series
derived in the previous section, the second-order results for the various
types of magic series are as follows:
\begin{eqnarray}
{\cal N}_{2}(N) &\approx& \left(\frac{1}{\pi} \sqrt{\frac{3}{e}} \right)\cdot \frac{\left(Ne\right)^{N}}{N^3} \left(1 + \frac{3}{5N} + 
\frac{31}{420N^2}\right).
\\
{\cal N}_{3}(N) &\approx& \left(\frac{\sqrt{3}}{\pi} \right) \cdot \frac{\left(N^2 e\right)^{N}}{N^4} 
\left(1 - \frac{11}{15N} + \frac{157}{126N^2}\right).
\\
{\cal N}_{4}(N) &\approx& \left(\frac{\sqrt{3}}{\pi} \right) \cdot \frac{\left(N^{3} e\right)^{N}}{N^{5}} 
\left(1 - \frac{7}{30N} - \frac{1249}{2520N^2}\right).
\\
{\cal N}_{\alpha > 4}(N) &\approx& \left(\frac{\sqrt{3}}{\pi} \right) \cdot \frac{\left(N^{\alpha-1} e\right)^{N}}{N^{\alpha + 1}} 
\left(1 - \frac{7}{30N} + \frac{11}{2520N^2}\right).
\end{eqnarray}
It is interesting to note that if the expansion of ${\cal N}_2(N)$ in decreasing powers of $N$
is moved to the denominator, one obtains a very slight correction to Bottomley's 
formula, Eq.~(\ref{eq:bottomley}), for the number of magic square series :
\begin{equation}
{\cal N}_{\rm{2}}(N) \approx \left(\frac{1}{\pi} \sqrt{\frac{3}{e}} \right)\cdot \frac{\left(Ne\right)^{N}}{N^3 - \frac{3}{5}N^2 + 
\left(\frac{2}{7} + \frac{1}{2100}\right)N}.
\end{equation}
The coefficient $(2/7)$ in Bottomley's empirical formula is replaced by $(2/7) + (1/2100)$.
All coefficients in the above formulae are exact;
all besides the $(3/5)$ in ${\cal N}_2(N)$ are novel, to the author's knowledge.

Table~\ref{table:magic-series} compares these second-order estimates with the exact counts, which are
known at least up to $N=1000$ for squares, $N=200$ for cubes, and $N=20$ for $4$-dimensional
hypercubes~\cite{trump-enum}.  The agreement is excellent: at these largest known values,
our second-order approximations are
correct to (respectively) $11$, $7$, and $4$ decimal places.  Moreover, the 
exact values allow us approximate the unknown third-order terms fairly well.  
Based on the contents of Table~\ref{table:magic-series}, the third-order terms are roughly 
$0.032/N^3$, $-1.44/N^3$, and $1.12/N^3$ for squares, cubes, and $4$-dimensional
hypercubes respectively.

\begin{table}
  \begin{center}
\begin{tabular}{|c|c|c|c|c|}
  \hline
   ${\alpha}$ & N & Est. vs. exact magic series & Rel. residual & Residual $\times N^3$ \\
  \hline
  2 & 500 & $1.148464537\underbar{053} \times 10^{1558}$ & $2.46\times 10^{-10}$ & 0.031 \\
     &        & $1.14846453733617811101 \times 10^{1558}$ & & \\
  \hline
  2 & 700 & $3.66527778\underbar{173} \times 10^{2286}$ & $9.00 \times 10^{-11}$ & 0.031 \\
     &  & $3.66527778205981116969 \times 10^{2286}$ & &\\
\hline
2 & 1000 & $6.591829225\underbar{191} \times 10^{3424}$ & $3.18 \times 10^{-11}$ & 0.032 \\
 &  & $6.59182922540146398832 \times 10^{3424}$ & & \\
  \hline
  \hline
  3 & 100 & $1.4713530\underbar{522}\times 10^{435}$ & $-1.43 \times 10^{-6}$ & -1.43 \\
    &  & $1.47135094282799112413\times 10^{435}$ & & \\
\hline
  3 & 150 & $1.01505\underbar{942}\times 10^{709}$ & $-4.25 \times 10^{-7}$ & -1.44 \\
    &  & $1.01505898472535873940 \times 10^{709}$ & & \\
\hline
  3 & 200 & $6.40624\underbar{667} \times 10^{997}$ & $-1.80 \times 10^{-7}$ & -1.44 \\
    &  & $6.40624551588569444014\times 10^{997}$ & & \\
  \hline
  \hline
  4 & 10 & $1.18\underbar{003}\times 10^{29}$ & $1.10\times 10^{-3}$ & 1.10\\
   &  & $1.18132052487666384802\times 10^{29}$ & & \\
 \hline
  4 & 15 & $1.95\underbar{748}\times 10^{53}$ & $3.28\times 10^{-3}$ & 1.11\\
    &  & $1.95811840766424991031\times 10^{53}$ & & \\
 \hline
   4 & 20 & $9.51\underbar{281}\times 10^{79}$ & $1.39\times 10^{-4}$ & 1.11 \\
   &  & $9.51413048876962267360\times 10^{79}$ & & \\
  \hline
\end{tabular}
\caption{{\bf Comparison of second-order magic series estimates to selected exact values.}  Results are
shown for the number of magic square series ($\alpha=2$), magic cube series ($\alpha=3$), and magic
series for $4$-dimensional hypercubes.  The underlined 
numbers are the first incorrect digits in each estimate.  The final column approximates the
coefficient of the $N^{-3}$ term for each $\alpha$; note that these approximations appear
to have converged.}
\label{table:magic-series}
\end{center}
\end{table}

\section{Extension to multimagic series}

\subsection{Bimagic series}

A \emph{bimagic square} is a magic square with the additional property that
the squares of the entries in each row, column, and diagonal have the same sum;
one then defines \emph{bimagic square series} as before.  Moreover, one can define
\emph{bimagic cube series}, etc., as well as trimagic and higher multimagic 
squares, cubes, etc., and their associated series.
The method described in this paper can be adapted to enumerate multimagic series
as well.  We limit ourselves here to briefly describing its application
to bimagic and trimagic series.

For bimagic series, we will want start by considering random vectors $(X_i, Y_i, Z_i)$ 
that are equal to $(1,i,\frac{1}{2}i(i-1))$ with probability $\beta$ and equal
to $(0,0,0)$ otherwise.  (If we were to choose $Z_i=i^2$ instead, $Y_i$ and $Z_i$ 
would be equal modulo $2$; in the $k$-space analysis, this would lead to multiple points of stationary 
phase, complicating the problem.)  The probability distribution of the sum of the first $m$
such vectors is related to the number of bimagic series just as in 
Eq.~(\ref{statform}).  This probability distribution approaches a 
three-dimensional Gaussian centered at
\begin{eqnarray}
\mu_x &=& \beta\sum_{i=1}^{m}1 = \beta m \\
\mu_y &=& \beta\sum_{i=1}^{m}i = \frac{1}{2}\beta m(m+1) \\
\mu_z &=& \beta\sum_{i=1}^{m}\frac{1}{2}i(i-1) = \frac{1}{6}\beta m(m^2-1).
\end{eqnarray}
Its covariance matrix has components
\begin{eqnarray}
\Sigma_{xx} &=& \beta(1-\beta) m \\
\Sigma_{xy} &=& \frac{1}{2}\beta(1-\beta)m(m+1) \\
\Sigma_{xz} &=& \frac{1}{6}\beta(1-\beta)m(m^2-1) \\
\Sigma_{yy} &=& \frac{1}{3}\beta(1-\beta)m(m+1)\left(m+\frac{1}{2}\right) \\
\Sigma_{yz} &=& \frac{1}{8}\beta(1-\beta)m(m^2-1)\left(m+\frac{2}{3}\right) \\
\Sigma_{zz} &=& \frac{1}{20}\beta(1-\beta)m(m^2-1)\left(m^2 - \frac{2}{3}\right)
\end{eqnarray}
To leading order in $1/m$, this matrix has determinant $\beta^{3} (1-\beta)^{3} m^9/8640$;
the Gaussian approximation for the peak probability is then
\begin{equation}
P^{*}_{g}(m,\beta) = \frac{1}{(2\pi)^{3/2}\sqrt{\det\Sigma}} \approx \frac{6\sqrt{30}}{\pi^{3/2}(1-\beta)^{3/2}\beta^{3/2}m^{9/2}}.
\end{equation}
For bimagic square series, the prefactor is the same as in 
Eq.~(\ref{prefactor-squares}), $\beta=1/N$, and $m=N^2$, so the
Gaussian approximation for the number of bimagic square series is
\begin{eqnarray}
{\cal N}_{2,g}^{(2)}(N) &=& \frac{(Ne)^N}{\sqrt{e}}\left(1 - \frac{1}{6N} + \cdots\right)\times \frac{6\sqrt{30}\left(1-\frac{1}{N}\right)^{-3/2}}{\pi^{3/2} N^{15/2}}  
\nonumber
\\
&=& \left(\frac{6}{\pi^{3/2}}\sqrt{\frac{30}{e}} \right)\cdot \frac{(Ne)^{N}}{N^{15/2}}\left(1 + \frac{4}{3N} +
\cdots\right).
\end{eqnarray}
For bimagic cube series, the prefactor is the same as in 
Eq.~(\ref{prefactor-cubes}), $\beta=1/N^2$, and $m=N^3$, so the Gaussian
approximation here is
\begin{eqnarray}
{\cal N}_{3,g}^{(2)}(N) &=& \left(N^2 e\right)^{N} \left(1 - \frac{1}{2N} + 
\cdots\right) \times \frac{6\sqrt{30}\left(1-\frac{1}{N^2}\right)^{-3/2}}{\pi^{3/2}N^{21/2}} 
\nonumber \\
&=& \left(\frac{6\sqrt{30}}{\pi^{3/2}} \right)\cdot \frac{\left(N^2 e\right)^{N}}{N^{21/2}}\left(1 - \frac{1}{2N} +
\cdots\right).
\end{eqnarray}

Perturbative corrections to the Gaussian results may be calculated exactly as before.  In 
fact, Eq.~(\ref{perturb-series}) is still correct, as is the expression
for the first correction coefficient $K_1$.  The only difference is that the 
propagator and vertex have changed.  The propagator is now
\begin{equation}
  \Pi^{(2)} = \left(
\begin{array}{ccc}
  9 & -36 & 60 \\
  -36 & 192 & -360 \\
  60 & -360 & 720 \\
\end{array}\right)
\end{equation}
up to corrections of order $1/m$.  The vertex is given by 
\begin{eqnarray}
\tilde{V}^{(2)}_{abc\ldots} &=& \frac{1}{m^{1+n_y(abc\ldots) + 2n_z(abc\ldots)}}\sum_{j=1}^{m} 
j^{n_y(abc\ldots)}\left(\frac{1}{2}j(j-1)\right)^{n_z(abc\ldots)} \nonumber \\
&\approx& \left(\frac{1}{2}\right)^{n_z(abc\ldots)}\frac{1}{1 + n_y(abc\ldots) + 
2n_z(abc\ldots)},
\end{eqnarray}
also to leading order in $1/m$. We find that
\begin{eqnarray}
K_1^{(2)} &=& -\frac{\tilde{V}_{abc}^{(2)} \tilde{V}_{def}^{(2)}}{8} \Pi_{ab}^{(2)}\Pi_{cd}^{(2)}\Pi_{ef}^{(2)} -\frac{\tilde{V}_{abc}^{(2)} \tilde{V}_{def}^{(2)}}{12}\Pi_{ad}^{(2)}\Pi_{be}^{(2)}\Pi_{cf}^{(2)}
+ \frac{\tilde{V}_{abcd}^{(2)}}{8} \Pi_{ab}^{(2)} \Pi_{cd}^{(2)} \nonumber
\\
&=& -\frac{1}{8}\cdot\frac{2781}{245} - \frac{1}{12}\cdot\frac{2403}{245} + \frac{1}{8}\cdot\frac{423}{35} 
\nonumber
\\ &=& -\frac{711}{980}.
\end{eqnarray}
Multiplying $1 + K_1^{(2)}/N$ by the Gaussian approximations gives
\begin{equation}
{\cal N}^{(2)}_{2}(N) \approx  \left(\frac{6}{\pi^{3/2}}\sqrt{\frac{30}{e}} \right)\cdot \frac{(Ne)^{N}}{N^{15/2}}\left(1 + 
\frac{1787}{2940N}\right)
\end{equation}
and
\begin{equation}
{\cal N}^{(2)}_{3}(N) \approx  \left(\frac{6\sqrt{30}}{\pi^{3/2}} \right)\cdot \frac{\left(N^2 e\right)^{N}}{N^{21/2}}\left(1 
- \frac{1201}{980N}\right)
\end{equation}
for the number of bimagic square series and bimagic cube series, respectively.

Table~\ref{table:bimagic-series} compares our approximations
to some known results (from \cite{boyer-enum}).
The number of bimagic square series is known to $N=28$, at which point our approximation
is correct to $4$ significant digits.  The number of bimagic cube series is only known to
$N=11$; our first-order approximation is no better than the Gaussian approximation at that 
value, but should improve thereafter.   We predict the number of bimagic square 
series of order $29$ to be $3.9714\times 10^{44}$, and the number of bimagic 
cube series of order $12$ to be $3.1966 \times 10^{20}$.

\begin{table}
  \begin{center}
\begin{tabular}{|c|c|c|c|c|}
  \hline
   $\alpha$ & N & Est. vs. exact bimagic series & Rel. Residual & Residual $\times N^2$ \\
\hline
   2 & 25 &  $7.68\underbar{395}\times 10^{35}$ & $9.45\times 10^{-5}$ & 0.059 \\
   &  & $7.68467875608077797721\times 10^{35}$ & & \\
  \hline
   2 & 26 &  $1.077\underbar{803}\times 10^{38}$ & $-6.50\times 10^{-4}$ & -0.44 \\
   &  & $1.07710220763567919575\times 10^{38}$ & & \\
  \hline
   2 & 27 &  $1.588\underbar{739}\times 10^{40}$ & $-4.35\times 10^{-4}$ & -0.32 \\
   &  & $1.58804753475269623163\times 10^{40}$ & & \\
  \hline
   2 & 28 &  $2.455\underbar{567}\times 10^{42}$ & $-4.49\times 10^{-5}$ & -0.035 \\
   &  & $2.45545658112397677916\times 10^{42}$ & & \\
  \hline
  \hline
  3 & 9 & $\underbar{5.93} \times 10^{11}$ & 0.098 & 7.98 \\
    &  & $6.51151145259 \times 10^{11}$ & & \\
\hline
  3 & 10 & $3.\underbar{607} \times 10^{14}$ & -0.038 & -3.75 \\
    &  & $3.47171191981324 \times 10^{14}$ & & \\
\hline
  3 & 11 & $\underbar{2.97}\times 10^{17}$ & 0.060 & 7.29 \\
    &  & $3.15035719463520007\times 10^{17}$ & & \\
  \hline
 \end{tabular}
\caption{{\bf Comparison of first-order bimagic series estimates to selected exact values.}  Results are
shown for the number of bimagic square series ($\alpha=2$) and bimagic cube series ($\alpha=3$).  The underlined 
numbers are the first incorrect digits in each estimate.  The final column approximates the
coefficient of the $N^{-2}$ term for each $\alpha$.}
\label{table:bimagic-series}
\end{center}
\end{table}

\subsection{Trimagic series}

A \emph{trimagic square} is a bimagic square with the additional property that
the cubes of the entries in each row, column, and diagonal have the same sum.
To enumerate trimagic series, we will consider random vectors $(X_i, Y_i, Z_i, W_i)$ 
that are equal to $(1,i,\frac{1}{2}i(i-1), \frac{1}{6}i(i-1)(i-2))$ with probability $\beta$ and equal
to $(0,0,0,0)$ otherwise.  All aspects of the preceding calculation go through 
as before, so we will omit the details.  One important difference is that the 
expected value of $\sum_{j=1}^{N^2}W_j$ is not an integer when 
$N=4k+2$, and so there are no trimagic series for $N$ of this form.  For all
other values of $N$, we find
\begin{equation}
{\cal N}_2^{(3)}(N) \approx \left(\frac{720}{\pi^2}\sqrt{\frac{105}{e}}\right)\cdot\frac{(Ne)^N}{N^{14}}
\end{equation}
to leading order.

The exact values (from \cite{boyer-enum}) are compared to the leading-order approximation in 
Table~\ref{table:trimagic-series}.  The agreement is fairly poor at the 
accessible values of $N$, and because the errors oscillate in magnitude, the $1/N$
correction cannot improve matters much.  There are clearly non-perturbative
effects that need to be better understood, and these effects are evidently more important
for higher multimagic series.  Indeed, similar (albeit smaller) 
oscillatory effects are already apparent in the bimagic series approximation.

\begin{table}
     \begin{center}
\begin{tabular}{|c|c|c|c|}
  \hline
   $\alpha$ & N & Est. vs. exact trimagic series & Exact/estimate \\
  \hline
  2 & 11 &  $2.04\times 10^{4}$  & $1.53$ \\
     &        & $31187$ & \\
  \hline
  2 & 12 &   $5.12\times 10^{5}$     & $4.35$ \\
     &        & $2226896$ & \\
  \hline
   2 & 13 &  $1.54 \times 10^{7}$ & $1.12$ \\
     &  &  $17265701$ & \\
\hline
   2 & 15 &  $2.22\times 10^{10}$ & $3.12$ \\
   &  & $69303997733$ & \\
  \hline
 \end{tabular}
\caption{{\bf Comparison of leading-order trimagic series estimates to selected exact values.}  Results are
shown for the number of trimagic square series ($\alpha=2$) only.}
\label{table:trimagic-series}
 \end{center}
\end{table}

\end{document}